\documentclass[graybox]{svmult}

\usepackage{mathptmx}       % selects Times Roman as basic font
\usepackage{helvet}         % selects Helvetica as sans-serif font
\usepackage{courier}        % selects Courier as typewriter font
\usepackage{type1cm}        % activate if the above 3 fonts are
\usepackage{makeidx}         % allows index generation
\usepackage{graphicx}        % standard LaTeX graphics tool
\usepackage{amssymb,amsmath}                                % when including figure files
\usepackage{multicol}        % used for the two-column index
\usepackage[bottom]{footmisc}% places footnotes at page bottom

\renewcommand{\le}{\leqslant}
\renewcommand{\ge}{\geqslant}
\renewcommand{\P}{\mathsf{P}}
\newcommand{\R}{\mathbb{R}}
\newcommand{\Rd}{\mathbb{R}^d}
\newcommand{\N}{\mathbb{N}}
\newcommand{\Z}{\mathbb{Z}}
\renewcommand{\E}{\mathsf{E}}
\newcommand{\ind}{{\bf{1}}}
\def\cov{\mathsf{c}\mathsf{o}\mathsf{v}}
\def\var{\mathsf{V}\mathsf{a}\mathsf{r}}
\newcommand{\Lip}{\mathsf{L}\mathsf{i}\mathsf{p}}
\newcommand{\dist}{\mathsf{d}\mathsf{i}\mathsf{s}\mathsf{t}}

\newcommand{\card}{\mbox{card}}
\newcommand{\tod}{\stackrel{d}{\longrightarrow}}

\makeindex             % used for the subject index
\begin{document}

\title*{Limit theorems for excursion sets of stationary random fields}
\titlerunning{LTs for excursion sets of random fields}
\author{Evgeny Spodarev}
\institute{Evgeny Spodarev \at Ulm University, Institute of
Stochastics, 89069 Ulm, Germany,
\email{evgeny.spodarev@uni-ulm.de}}

\maketitle

\abstract*{We give an overview of the recent asymptotic results on
the geometry of excursion sets of stationary random fields.
Namely, we cover a number of limit theorems of central type for
the volume of excursions of stationary (quasi--, positively or
negatively) associated  random fields with stochastically
continuous realizations for a fixed excursion level. This class
includes in particular Gaussian, Poisson shot noise, certain
infinitely divisible, $\alpha$--stable and max--stable random
fields satisfying some extra dependence conditions. Functional
limit theorems (with the excursion level being an argument of the
limiting Gaussian process) are reviewed as well. For stationary
isotropic $C^1$--smooth Gaussian random fields similar results are
available also for the surface area of the excursion set.
 Statistical tests of Gaussianity of a random field
which are of importance to real data analysis as well as results
for an increasing excursion level round up the paper.}

\abstract{We give an overview of the recent asymptotic results on
the geometry of excursion sets of stationary random fields.
Namely, we cover a number of limit theorems of central type for
the volume of excursions of stationary (quasi--, positively or
negatively) associated  random fields with stochastically
continuous realizations for a fixed excursion level. This class
includes in particular Gaussian, Poisson shot noise, certain
infinitely divisible, $\alpha$--stable and max--stable random
fields satisfying some extra dependence conditions. Functional
limit theorems (with the excursion level being an argument of the
limiting Gaussian process) are reviewed as well. For stationary
isotropic $C^1$--smooth Gaussian random fields similar results are
available also for the surface area of the excursion set.
Statistical tests of Gaussianity of a random field which are of
importance to real data analysis as well as results for an
increasing excursion level round up the paper.}

\section{Introduction}
\label{Spo:sec:intro}

Geometric characteristics such as Minkowski functionals (or
intrinsic volumes, curvature measures, etc.) of excursions of
random fields are widely used for data analysis purposes in
medicine (brain fMRI analysis, see e.g. \cite{ATW09},
\cite{SDT08}, \cite{TayWor07}, \cite{WorTay06}), physics and
cosmology (microwave background radiation analysis, see e.g.
\cite{MarPec11} and references therein), and materials science
(quantification of porous media, see e.g. \cite{MeckeStoyan02},
\cite{Torquato02}), to name just a few. Minkowski functionals
include the volume, the surface area and the Euler--Poincar\'e
characteristic (reflecting porosity) of a
 set with a sufficiently regular boundary.

Among the possible abundance of random field models, Gaussian
random fields are best studied  due to their analytic
tractability. A number of results starting with explicit
calculation of the moments of Minkowski functionals is available
for them since the mid seventies of the last century. We briefly
review these results in Section \ref{Spo:sec:MomentsVj}. However,
our attention is focused on the asymptotic arguments for (mainly
non--Gaussian) stationary random fields. There has been a recent
breakthrough in this domain starting with the paper
\cite{BuSpoTim12} where a central limit theorem (CLT) for the
volume of excursions of a large class of quasi--associated random
fields was proved. We also cover a number of hard--to--find
results from recent preprints and PhD theses.

The paper is organized as follows. After introducing some basic
facts on excursions and dependence structure of stationary random
fields in Section \ref{Spo:sec:prelim}, we briefly review the
limit theorems for excursions of  stationary Gaussian processes
($d=1$) in the next section. However, our focus is on the recent
results in the multidimensional case $d>1$ which is considered in
Sections \ref{Spo:sec:LTVol} and \ref{Spo:sec:LTSurf}. Thus,
Section \ref{Spo:sec:LTVol} gives (uni- and multivariate as well
as functional) central limit theorems for the volume of excursion
sets of  stationary (in general, non--Gaussian) random fields over
fixed, variable or increasing excursion levels. In Section
\ref{Spo:sec:LTSurf}, a similar scope of results is covered for
the surface area of the boundary of excursion sets of stationary
(but possibly anisotropic)  Gaussian random fields in different
functional spaces. The paper concludes with a number of open
problems.

\section{Preliminaries}
\label{Spo:sec:prelim}

Fix a probability space $(\Omega, \cal F, \P)$. Let $X=\{
X(t,\omega), \, t \in \Rd,\, \omega\in\Omega \}$ be a stationary
(in the strict sense) real valued measurable (in $(t,\omega)\in
\Rd\times \Omega$) random field. Later on we suppress $\omega$ in
the notation. For integrable $X$ we assume $X$ to be centered
(i.e., $\E X(o)=0$ where $o\in\Rd$ is the origin point). If the
second moment of $X(o)$ exists then we denote by $C(t)=\E \,
\left( X(o) X(t) \right)$, $t\in\Rd$ the covariance function of
$X$.

Let  $ \| \cdot\|_2$ be the Euclidean norm  in $\Rd$ and $\dist_2$
the Euclidean distance: for two sets $A,B\subset \Rd$ we put
$\dist_2(A,B)=\inf \{ \| x-y\|_2:\, x\in A, y\in B\}$. Denote by $
\| \cdot\|_\infty $ the supremum norm in $\Rd$ and by
$\dist_\infty$ the corresponding distance function.

Let $\tod$ mean convergence in distribution. Denote by $A^c$ the
complement and by $\mbox{int}(A)$ the interior of a set $A$ in the
corresponding ambient space which will be clear from the context.
Let $\card (A)$ be the cardinality of a finite set $A$. Denote by
$B_r(x)$ the closed Euclidean ball with center in $x\in\Rd$ and
radius $r>0$. Let ${\cal H}^{k}(\cdot)$ be the $k$--dimensional
Hausdorff measure in $\Rd$, $0\le k\le d$. In the sequel, we use
the notation $\kappa_{j}={\cal H}^{j}(B_1(o))$, $j=0,\ldots, d$.

To state limit theorems, one has to specify the way of expansion
of windows $W_n\subset T$, where the random field $X=\{ X(t), \;
t\in T \}$ is observed, to the whole index space $T=\Rd$ or
$\Z^d$. A sequence of compact Borel sets $ (W_n)_{n\in\mathbb{N}}$
is called a \emph{Van Hove sequence} (\emph{VH})\index{Van Hove
sequence} if $W_n\uparrow\R^d$ with
$$ \lim_{n\rightarrow\infty}V_d \left(W_n\right) = \infty\ \ \text{  and  }\ \ \lim_{n\rightarrow\infty}\frac{V_d \left(\partial W_n \oplus B_r(o)\right)}{V_d \left(W_n\right)} = 0,\ \ r>0.$$
A sequence of finite subsets $U_n\subset \Z^d$, $n\in\N$ is called
\emph{regular growing}\index{regular growth} if
$$\card(U_n)\to\infty \quad \mbox{and}\quad \card(\delta U_n)/\card(U_n)\to 0  \; \mbox{
as } \; n\to\infty$$ where $\delta U_n=\{ j\in\Z^d\setminus U_n:
\; \dist_\infty (j,U_n)=1 \}$ is the discrete boundary of $U_n$ in
$\Z^d$.

\subsection{Excursion sets and their intrinsic volumes} \label{Spo:sec:prelim:subsec:excurs}

The \emph{excursion set}\index{excursion set} of $X$ at level
$u\in\R$ in the compact observation window $W\subset \Rd$ is given
by $A_u(X,W)=\{t\in W: X(t)\ge u   \}$. The \emph{sojourn
set}\index{sojourn set} under the level $u$ is $S_u(X,W)=\{t\in W:
X(t)\le u   \}$, respectively.

 Due to measurability of $X$, $A_u(X,W)$ and $S_u(X,W)$ are
random Borel sets. If $X$ is a.s. upper (lower) semicontinuous
then $A_u(X,W)$ ($S_u(X,W)$, respectively) is a random closed set
(cf. \cite[Section 5.2.1]{Molch05}).

A popular way to describe the geometry of excursion sets is via
their \emph{intrinsic volumes}\index{intrinsic volumes} $V_j$,
$j=0,\ldots,d$. They can be introduced for various families of
sets such as convex and polyconvex sets \cite[Chapter 4]{schn93},
sets of positive reach and their finite unions \cite{fed59},
unions of basic complexes \cite[Chapter 6]{AdlerTaylor07}. One
possibility to define $V_j(K)$, $j=0,\ldots,d$ for a set $K$
belonging to the corresponding family is given by the
\emph{Steiner formula}  (see e.g. \cite[Section 13.3]{sant76}) as
the coefficients in the polynomial expansion of the volume of the
tubular neighbourhood $K_r=\{x\in\Rd: \dist_2 (x,K)\le r\}$ of $K$
with respect to the radius $r>0$ of this neighbourhood:
$$
{\cal H}^{d}\left(K_r\right)=\sum_{j=0}^d \kappa_{d-j} V_{j}(K)
r^{d-j}
$$
for admissible $r>0$ (for convex $K$, these are all positive $r$).
 The
geometric interpretation of intrinsic volumes $V_j(K)$,
$j=1,\ldots, d-2$ can be given in terms of integrals of elementary
symmetric polynomials of principal curvatures for convex sets $K$
with $C^2$--smooth boundary, cf. \cite[Sections 13.5-6]{sant76}.
Without going into details here, let us discuss the meaning of
some of $V_j\left(A_{u}(X,W)\right)$, $j=0,\ldots,d$ in several
dimensions.

For $d=1$, $V_{1} \left(A_{u}(X,W)\right) $ is the length of
excursion intervals and  $V_{0} \left(A_{u}(X,W)\right) $ is the
number of upcrossings of level $u$ by the random process $X$
within $W$.

For dimensions $d\ge 2$, $V_d \left(A_{u}(X,W)\right)$ is always
the volume (i.e., the Lebesgue measure) of $A_{u}(X,W)$ and
$V_{d-1} \left(A_{u}(X,W)\right) $ is half the surface area, i.e,
$1/2 \cdot {\cal H}^{d-1}\left(\partial A_{u}(X,W)\right)$. The
\emph{Euler characteristic} $V_{0} \left(A_{u}(X,W)\right)$ is a
topological measure of ``porosity'' of excursion set $A_{u}(X,W)$.
For ``basic'' sets $A$ (e.g., non--empty convex sets or sets of
positive reach) we set $V_{0}(A)=1$. Then $V_0$ is defined for
unions of basic sets by additivity. One can show that for $d=2$ it
holds
$$V_{0}(A)= \card \{ \mbox{connented components of }A
  \}- \card \{ \mbox{holes of }A\}.$$

The existence of $V_{j} \left(A_{u}(X,W)\right)$, $j=d,d-1$, is
clear since $A_{u}(X,W)$ is a Borel set whose Lebesgue and
Hausdorff measures are well defined. Intrinsic volumes $V_j$ of
lower orders $j=0,\ldots, d-2$ are well defined e.g. for excursion
sets of sufficiently smooth (at least $C^2$) deterministic
functions (cf. \cite[Theorem 6.2.2]{AdlerTaylor07} ) and Gaussian
random fields (cf. \cite[Theorem 11.3.3]{AdlerTaylor07})
satisfying some additional conditions.

\subsection{Dependence concepts for random fields} \label{Spo:sec:prelim:subsec:depend}

To prove limit theorems for a random field $X$, some conditions
have to be imposed on the structure of the dependence of $X$.
Mixing conditions that are usually required (cf. e.g.
\cite{Doukhan94}, \cite{Bra07}) are however rather difficult to
check for a particular random field under consideration. For this
practical reason, we follow the books \cite{Bulinski},
\cite[Chapter 10]{spodLNM} and introduce \emph{association} as
well as related dependence concepts.

A random field $X=\left\{X(t),\; t\in \Rd\right\}$  is called
\emph{ associated ({\bf A})}\index{association} if
$$
 \cov\left(f\left(X_I\right),g\left(X_I\right)\right) \ge 0
$$
for any finite  subset $I\subset \Rd$, and for any bounded
coordinatewise non--decreasing functions $f:\R^{\card
(I)}\rightarrow\R$, $g:\R^{\card(I)}\rightarrow\R$ where $X_I=\{
X(t),\; t\in I \}$.

 A random field $X=\left\{X(t),\; t\in \Rd\right\}$  is called
 \emph{ positively (${\bf PA}$)} or \emph{negatively (${\bf NA}$) associated}\index{positive, negative association} if
$$
 \cov\left(f\left(X_I\right),g\left(X_J\right)\right) \ge 0 \quad
 (\le 0, \; \mbox{resp.})
$$
for all finite disjoint subsets $I,J\subset \Rd$, and for any
bounded coordinatewise non--decreasing functions $f:\R^{\card
(I)}\rightarrow\R$, $g:\R^{\card(J)}\rightarrow\R$. It is clear
that if $X\in{\bf A}$ then $X\in{\bf PA}$.

Subclasses of ${\bf A}$ (${\bf PA},$ ${\bf NA}$)-- fields are
certain infinitely divisible (e.g., max-stable and
$\alpha$-stable) random fields. In particular, a Gaussian random
field with non--negative covariance function is associated.

A random field $X=\left\{X(t),t\in \Rd\right\}$ with finite second
moments is called \emph{quasi-associated}\index{quasi-associated
random field} (${\bf QA}$) if
$$
 \left|\cov\left(f\left(X_I\right),g\left(X_J\right)\right)\right| \leq \sum_{i\in I} \sum_{j\in J}\Lip_i \left(f\right)\Lip_j\left(g\right) \left|\cov\left(X\left(i\right),X\left(j\right)\right)\right|
$$
for all finite disjoint subsets $I,J\subset \Rd$, and for any
Lipschitz functions $f:\R^{\card (I)}\rightarrow\R$,
$g:\R^{\card(J)}\rightarrow\R$ where $\Lip_i\left(f\right)$ is the
Lipschitz constant of function $f$ for coordinate $i$. It is known
that if square integrable $X\in {\bf A} ({\bf PA}, {\bf NA})$ then
$X\in {\bf QA}$, cf. \cite[Theorem 5.3]{Bulinski}.

A real-valued random field $ X=\{X\left(t\right),t\in \R^d\} $ is
called \emph{$ \left(BL,\theta\right) $-dependent}\index{$
\left(BL,\theta\right) $-dependence} if there exists a
non--increasing sequence $ \theta = \{\theta_r\}_{r\in\R_0^+}$,
$\theta_r \downarrow 0$ as $r\rightarrow\infty$ such that for any
finite disjoint sets $ I, \ J \subset \Rd $ with $
\dist_\infty\left(I,J\right) = r\in\R_0^+ $, and any bounded
Lipschitz functions $f:\R^{\card (I)}\rightarrow\R$,
$g:\R^{\card(J)}\rightarrow\R$, one has
$$ \left|\cov\left(f\left(X_I\right),g\left(X_J\right)\right)\right| \leq \sum_{i\in I} \sum_{j\in J}\Lip_i \left(f\right)\Lip_j\left(g\right) \left|\cov\left(X\left(i\right),X\left(j\right)\right)\right|\theta_r. $$
It is often possible to choose $\theta$ as the \emph{Cox--Grimmett
coefficient}\index{Cox--Grimmett coefficient}
$$ \theta_r = \sup_{y\in\R^d}\int_{\R^d \setminus B_{r}^\infty (y)} \left|
\cov\left(X\left(y\right),X\left(t\right)\right)\right|dt $$ where
$B_{r}^\infty (y)=\{ x\in\Rd: \; \| x-y\|_\infty \le r\}$. It can
be easily seen that if $X\in {\bf QA}$ and its covariance function
is absolutely integrable on $\Rd$ then $X$ is
$\left(BL,\theta\right)$--dependent.

\section{Excursions of stationary Gaussian processes}
\label{Spo:sec:ExcStochProc}

Excursions of stochastic processes is a popular research topic in
probability theory since many years, see e.g. \cite{Berman92}
% \cite{DOT03}
and references in \cite{ILRS13}. The vast literature on this
subject for different classes of processes such as L\'evy,
diffusion, stable, Gaussian ones, etc. can be hardly covered by
one review. For this reason, we concentrate on the excursions of
(mainly stationary) Gaussian processes here.

Let $X=\{ X(t), \, t\ge 0 \}$ be a centered real valued Gaussian
process. If $X$ is a polynomial of degree $n$ with iid
$N(0,1)$-distributed coefficients then the mean number of real
roots of the equation $X(t)=0$ was first obtained by M. Kac
\cite{Kac43}. It initiated a substantial amount of papers on the
roots of random algebraic polynomials, see \cite{BordCha12} for a
 review. For $C^1$--smooth stationary
Gaussian processes $X$, expectation of the number of upcrossings
of a level $u$ by $X$ in time interval $[0,1]$ has been studied in
\cite{Rice44,Rice45}, \cite{Bul61}, etc. Higher order factorial
moments are considered in \cite{CramerLead65}, see also references
therein, and \cite{Bel66,Bel67}. For reviews (also including
results on non--Gaussian stationary processes) see \cite[Sections
7.2-7.3]{LeLiRoo83} and \cite[Chapter 3]{Azais}. In \cite{Adler76}
 and  \cite{Adler81}, the notion
of the number of upcrossings\index{number of upcrossings} of level
$u$ for random processes has been generalized to the
Euler-Poincar\'e characteristic of excursion sets of random
fields.

The first proof of a central limit theorem for the number
 of zeros of a stationary Gaussian process within an increasing time interval was given in
\cite{Malev69}. Cuzick \cite{Cuzick76} refined the assumptions
given in
 \cite{Malev69} and proved a central limit theorem for the number
 of zeros $N_X(T)=2 V_0 (A_0(X; [0,T]))$ of a centered separable stationary Gaussian process $X=\{ X(t), \, t\ge 0 \}$ in the time interval $[0,T]$ as well as
 analogous results for integrals $\int_0^T g\big(X(t)\big)\, dt$
 as $T\to\infty$. He used approximations by $m$--dependent random
 processes with spectral representation as a method borrowed from \cite{Malev69}.
In more detail, let $C(t)$ be twice differentiable with $C(0)=1$,
$C^{\prime\prime}(0)=-\lambda_2$ and variogram $\gamma$ of
$X^\prime$ be given by
$\gamma(h)=C^{\prime\prime}(h)-C^{\prime\prime}(0)=1/2 \E
(X^\prime(h)-X^\prime(0))$, $h\ge 0$.
\begin{theorem}[\cite{Cuzick76}]
If $C$, $C^{\prime\prime}$ are square integrable on $\R_+$,
$\int_0^\varepsilon \gamma(t)/t \, dt<\infty$ for some
$\varepsilon>0$  and
\begin{equation}\label{Spo:eq:VarCond}
\var N_X(T)/ T \to \sigma_2>0 \; \mbox{ as } \;  T\to + \infty
\end{equation} then
$$
T^{-1/2}\left( N_X(T) -\E\, N_X(T)\right) \tod N(0,\sigma^2) \;
\mbox{ as } \;  T\to + \infty
$$
where $$\sigma^2=\pi^{-1} \left( \lambda_2^{1/2}+
\int\limits_0^\infty \left( \frac{ \E\,\left(
|X^{\prime}(0)X^{\prime}(t)| | X(0)=X(t)=0\right)
}{\sqrt{1-C^2(t)}}- (\E\, |X^{\prime}(0)|)^2 \right)\,dt
\right).$$
\end{theorem}
Condition \eqref{Spo:eq:VarCond} is difficult to check and is
substituted in \cite[Lemma 5]{Cuzick76} by a more tractable
sufficient condition involving $C$ and $\lambda_2$. Piterbarg
\cite{Piterbarg78} managed to prove the above theorem by
substituting condition \eqref{Spo:eq:VarCond} with
$$
\int\limits_0^\infty t \left( |C(t)|+|C^\prime (t)|+|C^{\prime
\prime}(t)| \right)\,dt<\infty.
$$
He approximates the point process of upcrossings of $X$ of level
$u$ by a strongly mixing point process.

\begin{theorem}[\cite{Cuzick76}]
Let $X$ be a stationary Gaussian process with covariance function
 $C$ being integrable on $\R_+$. For any measurable function $g:
 \R\to\R$ such that $\E\, g^2(X(0))<\infty$ and $g(x)-g(0)$ is not
  odd it holds
\begin{equation}\label{eq:CLTFktCuzick}
 T^{-1/2}\left( \int_0^T g\big(X(t)\big)\,
dt - T \E\, g\big(X(0)\big)\right) \tod N(0,\sigma^2) \; \mbox{ as
} \; T\to + \infty
\end{equation}
where $\sigma^2>0$.
\end{theorem}
It is clear that the choice $g(x)=\ind \{x\in\R: x\ge u  \}$ for
any $u\in\R$ leads to the central limit theorem for the length
$V_1 (A_u(X; [0,T]))$ of excursion intervals of $X$ in $[0,T]$.

Elizarov \cite{Eliz88} first proved a functional central limit
theorem for the sojourn times of the stationary Gaussian process
under the level $u$, in our terms, for $V_1(S_u(X; [0,T]))$  if
excursion level $u$ is allowed to vary within $\R$. Additionally,
an analogous result for local times
$$
\lim_{\varepsilon\to+0}\frac{1}{2\varepsilon}\left(
V_1(S_{u+\varepsilon}(X; [0,T])) - V_1(S_{u-\varepsilon}(X;
[0,T])) \right)
$$
was given. Both results were proved in the functional space
$C[0,1]$ after the substitution $u\mapsto f(x)$, $x\in [0,1]$
where $f\in C[0,1]$ is a monotonously increasing function with
$f(0)=-\infty$, $f(1)=\infty$.

Belyaev and Nosko \cite{Bel69} proved limit theorems for
$V_1(A_u(X; [0,T]))$, $T\to\infty$ as $u\to\infty$ for stationary
ergodic processes $X$ satisfying a number of additional (quite
technical) assumptions. In particular, these assumptions are
satisfied if $X$ is an ergodic Gaussian stationary process with
twice continuously differentiable covariance function such that
$$\left| C^{\prime \prime}(t)-
C^{\prime \prime}(0)  \right| \le a/|\log |t||^{1+\varepsilon},
\quad t\downarrow 0$$ for some constants $a, \varepsilon>0$.

Slud \cite{Slud94} gave a multiple Wiener- It\^{o} representation
for the number of crossings of  a $C^1$--function $\psi$ by $X$.
In  \cite{KratzLeon01}, methods of  \cite{Malev69} and
\cite{Cuzick76} are generalized to the case of functionals of $X$,
$X^\prime$ and  $X^{\prime\prime}$. CLTs for the number
 of  crossings of  a smooth curve $\psi$ by a Gaussian process $X$  as well as for the number of specular points of $X$ (if $X$ is a Gaussian process in time and space) are given in \cite{KratzLeon10}.
For a  review of results on moments and limit theorems for
different characteristics of stationary Gaussian processes see
\cite{Kratz06}. In \cite{ILRS13},  CLTs for the multivariate
non--linear weighted functionals (similar to those in
\eqref{eq:CLTFktCuzick}) of Gaussian stationary processes
 with multiple singularities in their spectra, having a
covariance function  belonging to a certain  parametric family,
are proved.

\section{Moments of $V_j \left(A_u(X,
W) \right)$ for Gaussian random fields} \label{Spo:sec:MomentsVj}

We briefly review the state of the art for $\E\, V_j \left(A_u(X,
W) \right)$ of Gaussian random fields $X$. For recent extended
surveys see the books \cite{AdlerTaylor07} and  \cite{Azais}. For
stationary (isotropic) Gaussian fields $X$, stratified
$C^2$--smooth compact manifolds $W\subset \Rd$ and any $u\in\R$,
formulae for  $\E \, V_j \left(A_u(X, W) \right)$, $j=0,\ldots, d$
are given in \cite[Theorems 13.2.1 and 13.4.1]{AdlerTaylor07}.

Apart from obtaining exact (or asymptotic as $u\to\infty$)
formulae for $\E \, V_j \left(A_u(X, W) \right)$, $j=0,\ldots, d$,
the possibility of an estimate
\begin{equation}\label{eq:EPEur}
 \left| \P\left(\sup_{t\in W} X(t)>u \right)- \E\, V_0
 \left(A_u(X, W)
\right)  \right|\le g(u)
\end{equation}
(the so--called \emph{Euler-Poincar\'e
heuristic}\index{Euler-Poincar\'e heuristic}) with $g(u)=o(1)$ as
$u\to\infty$  is of special interest. It has been proved in
\cite[Theorem 14.3.3]{AdlerTaylor07} with $g(u)=c_0 \exp\{-u^2
(1+\alpha)/2\}$ for some positive constants $c_0$ and $\alpha$ if
$X$ is a (non)stationary Gaussian random field with constant
variance on a stratified manifold $W$ as $u\to\infty$.
 Lower and upper bounds for the density of
supremum of stationary Gaussian random fields $X$ (which imply
relation \eqref{eq:EPEur}) for any $u\in\R$ are given in
\cite[Theorem 8.4]{Azais}. Similar bounds are  proven in
\cite[Theorem 8.10]{Azais} for non--stationary Gaussian random
fields $X$ with a unique point of maximum of variance in
$\mbox{int}(W)$ as $u\to\infty$.

In \cite{SpoZap12}, asymptotic behavior of $\E\, V_j \left(A_u (X,
[a,b]^d) \right)$, $j=0,d-1,d$ of non--stationary sufficiently
smooth Gaussian random fields is studied as the excursion level
$u\to\infty$. The variance of these fields is assumed to attain a
global maximum at a vertex of $[a,b]^d$. It is shown that the
heuristic \eqref{eq:EPEur} still holds true.

A interesting rather general formula for the mean surface area of
Gaussian excursions \index{mean surface area of Gaussian
excursions} is proven in \cite{IbrZap10}.
 Let $W$ be a compact subset of $\Rd$ with a
non-empty interior and a finite Hausdorff measure of the boundary.
Let $X=\{X(t), \; t\in W\}$ be a Gaussian random field  with mean
$\mu(t)=\E X(t)$ and variance $\sigma^2(t)=\var X(t)$. For an
arbitrary (but fixed) excursion level $u\in\R$ introduce the zero
set $ \nabla_{X}^{-1}(0) $ of the gradient of the normalized field
$(X-u)/\sigma$ by $ \nabla_{X}^{-1}(0)=\{t\in W\,:\, \nabla\left(
(X(t)-u)/\sigma(t)\right)=0\}.$
\begin{theorem}[\cite{IbrZap10}]\label{1822}
Assume that  $X\in {C}^1(W)$ a.s.,  $\E V_{d-1}\left(
\nabla_{X}^{-1}(0)  \right)<\infty$ and $\sigma(t)>0$ for all
$t\in W$. Then
$$
\E V_{d-1}\left(
\partial A_u(X,W)  \right)=\frac{1}{2\sqrt{2\pi}}\int_W\,\exp\left[-\frac{(\mu(t)-u)^2}{2\sigma^2(t)}\right]\E\bigg\|\nabla
\left( (X(t)-u)/\sigma(t)\right) \bigg\|_2 \,dt.
$$
\end{theorem}

Asymptotic formulae for $\E \, V_j \left(A_u(X, W) \right)$,
$j=0,\ldots, d$ as $u\to\infty$  of three subclasses of stable
random fields (subgaussian, harmonizable,
concatenated--harmonizable ones) are given in \cite{AST10}.

\section{Volume of excursion sets of stationary random fields}
\label{Spo:sec:LTVol}\index{volume of excursion sets}

The first limit theorems of central type for the volume of
excursion sets (over a fixed level $u$) of stationary isotropic
Gaussian random fields were proved in \cite[Chapter 2]{IvLeon89}.
There, the case of short and long range  dependence (Theorem 2.2.4
and Example 2.2.1, Theorem 2.4.6) was considered. The CLT followed
from a general Berry-Ess\'een-type bound for the distribution
function of properly normed integral functionals
\begin{equation}\label{eq:IntFkt}
\int_{B_r(o)} G\big(X(t)\big) \, dt
\end{equation}
as $r\to\infty$ where $G:\R\to\R$ is a function such that
$\E\,G^2(X(o))<\infty$ satisfying some additional assumptions, cf.
also \cite{Leon88}. To get the volume $V_d \left( A_{u}\left(X,
B_r(o)\right) \right)$ out of \eqref{eq:IntFkt}, set
$G(x)=\ind(x\ge u)$. The isotropy of $X$ was essential as one used
expansions with respect to the basis of Chebyshev-Hermite
polynomials in the proofs. The cases of
 $$G(x)=\ind(|x|\ge u),\; \max\{0,x \}, \; |x|$$ as well as of $G$
 depending on a parameter and of weighted integrals in
 \eqref{eq:IntFkt} are considered as well.

In a remark \cite[p. 81]{IvLeon89}, it was noticed that similar
CLTs can be expected for non-Gaussian mixing random fields. The
aim of this Section is to review the recent advances in proving
such CLTs for various classes of stationary random fields that
include also the (not necessarily isotropic) Gaussian  case.

For instance,  random fields with singularities of their spectral
densities are considered in \cite{Leon99}. In Section 3.2 of that
book, non--central limit theorems for the volume of excursions of
stationary isotropic Gamma correlated and $\chi^2$-random fields
over a radial surface (i.e., the level $u$ is not constant
anymore, but a function of $\|t\|_2$, where $t\in\Rd$ is the
integration variable in \eqref{eq:IntFkt}) are proved.
(Non)central limit theorems for functionals \eqref{eq:IntFkt} of
stationary isotropic vector--valued Gaussian random fields are
given in the recent preprint \cite{LeonOl13}. There, the case of
long and short range dependence is considered as well as
applications to $F$-- and $t$--distributed random fields.

The asymptotic behavior of tail probabilities
$$
\P \left( \int_{W} e^{X(t)} \, dt>x  \right), \quad x\to\infty
$$
for a homogeneous smooth Gaussian random field $X$ on a compact
$W\subset \Rd$ is considered in \cite{Liu12}, see \cite{LiuXu12}
for further extensions.

\subsection{Limit theorems for a fixed excursion level}
\label{Spo:sec:LTVol:subsec:fixedLevel}

The main result (which we call a \emph{methatheorem}) can be
formulated as follows.
\begin{theorem}[Methatheorem]\label{theo:methaCLT}
 Let $X$ be a strictly stationary random field satisfying
\emph{some additional conditions} and $u\in\R$ fixed. Then, for
any sequence of $VH$-growing sets $W_n \subset \R^d $, one has
\begin{equation}\label{eq:methaCLT}
    \frac{V_d \left( A_{u}\left(X,W_n\right)\right) - \P(X(o)\geq u)\cdot V_d \left(W_n\right)}{\sqrt{V_d \left(W_n\right)}} \xrightarrow{\mathsf{d}} \mathcal{N}\left(0,\sigma^2(u)\right)
\end{equation}
as $ n \rightarrow \infty $. Here
\begin{equation}\label{eq:methaCLTassymptVar}
   \sigma^2(u) =\int_{\R^d} \cov\left(\ind\{X\left(o\right)\geq u\},\ind\{X\left(t\right)\geq u\}\right)\,dt.
\end{equation}
\end{theorem}
Depending on the class of random fields, these additional
conditions will vary. First we consider the family of
square--integrable random fields.

\subsubsection{Quasi-associated random fields} \label{Spo:sec:LTVol:subsect:fixedLevel:finSecMom}

\begin{theorem}[\cite{BuSpoTim12}]\label{th:QACLT}
Let $ X=\{X\left(t\right),\, t\in\R^d\} \in {\bf QA}$ be a
stationary square-integrable random field with a continuous
covariance function $C$ such that $ |C(t)| =
\mathcal{O}\left(\left\|t\right\|_2^{-\alpha}\right)$ for some
$\alpha>3d$ as  $\left\|t\right\|_2\rightarrow\infty$. Let $X(o)$
have a bounded density. Then $\sigma^2(u)\in(0,\infty)$ and
Theorem \ref{theo:methaCLT} holds true.
\end{theorem}

Let us give an idea of the proof. Introduce the random field $Z=\{
Z(j),\; j\in\Z^d \}$ by
\begin{equation}\label{eq:FieldZ}
Z(j)=\int_{j+[0,1]^d} \ind \big\{ X(t)\ge u\big\} \, dt - \Psi(u),
\quad j\in\Z^d.
\end{equation}
Here $ \Psi(u)=\P \big( X(o)> u\big) $ is the tail distribution
function of $X(o)$. It is clear that the sum of $ Z(j)$ over
indices $j\in W_n\cap \Z^d $ approximates the numerator in
\eqref{eq:methaCLT}. One has to show that $Z$ can be approximated
by a sequence of  $ \left(BL,\theta\right) $-dependent stationary
centered square-integrable random  fields $Z_\gamma$,
$\gamma\downarrow 0$, on $\Z^d$. The proof finishes by applying
the following CLT to $Z_\gamma$ for each $\gamma>0$.

 \begin{theorem}[\cite{Bulinski}, Theorem 3.1.12]\label{theo:BLThetaCLT}
  Let $ Z=\{Z(j),j\in\Z^d\} $ be a $ \left(BL,\theta\right) $-dependent strictly stationary centered square-integrable random field. Then, for any sequence of regularly growing sets $ U_n \subset \Z^d $, one
  has
$$ S\left(U_n\right)/\sqrt{\card\left(U_n\right)} \xrightarrow{\mathsf{d}} \mathcal{N}\left(0,\sigma^2\right) $$
as $ n \rightarrow \infty $, with
$$
   \sigma^2 = \sum_{j\in\Z^d} \cov\left(Z\left(o\right),Z\left(j\right)\right).
$$
 \end{theorem}

We give two examples of random fields satisfying Theorem
\ref{th:QACLT}.
\begin{example}[\cite{BuSpoTim12}]\label{ex:ShotNoise}
  Let $ X=\{X(t),\; t\in\R^d\}$ be a stationary \emph{shot noise random field}\index{shot noise random field}  given by
    $X(t)=\sum_{i\in\N} \xi_i \varphi(t-x_i)$ where
   $\Pi_\lambda=\{x_i\}$ is a
homogeneous Poisson point process in $\R^d$ with intensity
$\lambda\in(0,\infty)$, $\{\xi_i\}$ is a family of i.i.d.
non--negative random variables with $\E \, \xi_i^2<\infty$ and
characteristic function $\varphi_\xi$. Assume that $\Pi_\lambda$
and $\{\xi_i\}$ are independent. Moreover, let
$\varphi:\R^d\to\R_+$ be a bounded and uniformly continuous Borel
function with   $\varphi(t)\leq g_0(\|t\|_2) =
\mathcal{O}\left(\|t\|_2^{-\alpha}\right)$ as $\|t\|_2\to\infty$
for a function $g_0:\R_+\to\R_+$, $\alpha>3d$, and
$$
\int\limits_{\R^d} \left|\, \exp\left\{\lambda \int_{\R^d}
\left(\varphi_\xi(s\varphi(t))-1\right)\,dt\right\} \right|\,
ds<\infty.
$$
Then Theorem \ref{th:QACLT} holds true.
\end{example}

\begin{example}[\cite{BuSpoTim12}]\label{ex:Gauss}
   Consider a stationary Gaussian random field \index{Gaussian random field} $ X=\{X\left(t\right),t\in\R^d\}$
   with a continuous covariance function $C(\cdot)$ such that
   $ |C(t)| = \mathcal{O}\left(\left\|t\right\|_2^{-\alpha}\right)$ for some $\alpha>d$ as
   $\left\|t\right\|_2\rightarrow\infty$.
   Let $ X\left(o\right)\sim \mathcal{N}\left(a,\tau^2\right)$.
   Then,
Theorem \ref{th:QACLT} holds true with
$$
   \sigma^2(u) = \frac{1}{2\pi}\int_{\R^d} \int_0^{\rho(t)} \frac{1}{\sqrt{1-s^2}}\,e^{-\frac{(u-a)^2}{\tau^2\left(1+s\right)}}\,ds\,dt,
$$
where $\rho(t) =
\mathsf{c}\mathsf{o}\mathsf{r}\mathsf{r}(X(o),X(t)).$ In
particular, for $u=a$ one has
$$\sigma^2(a) = \frac{1}{2\pi}\int_{\R^d} \arcsin\left(\rho(t)\right)\,dt. $$
\end{example}

\subsubsection{PA- or NA-random fields} \label{Spo:sec:LTVol:subsect:fixedLevel:Stoch Cont}

What happens if the field $X$ does not have the finite second
moment? In this case, another set of conditions for our
methatheorem to hold was proven in \cite[Theorem 3.59]{Karcher12}.

\begin{theorem}\label{theo:PANA_CLT}
 Let $ X=\{X(t),\; t\in\R^d\} \in {\bf PA} ({\bf NA})$ be stochastically continuous satisfying the following properties:
  \begin{enumerate}
     \item the asymptotic variance  $\sigma^2(u)\in(0,\infty)$ (cf. its definition in \eqref{eq:methaCLTassymptVar}),
     \item $\P \left(X(o)=u \right)=0$  for the chosen level
     $u\in\R$.
   \end{enumerate}
   Then Theorem \ref{theo:methaCLT} holds.
\end{theorem}
The idea of the proof is first to show that the random field $Z=\{
Z(j),\; j\in\Z^d \}$ defined in \eqref{eq:FieldZ}
 is {\bf PA} ({\bf NA}).  Second,
use \cite[Theorem 1.5.17]{Bulinski} to prove that $Z$ is
$\left(BL,\theta\right)$--dependent. Then apply Theorem
\ref{theo:BLThetaCLT} to $Z$.

A number of important classes of random fields satisfy Theorem
\ref{theo:PANA_CLT}. For instance, stationary infinitely divisible
random fields $X=\{ X(t), \; t\in\R^d\}$ with spectral
representation
$$
X(t)=\int_E f_t(x)\, \Lambda (dx), \quad t\in\Rd,
$$
where $\Lambda$ is a centered independently scattered infinitely
divisible random measure on space $E$ and $f_t:E\to\R_+$ are
$\Lambda$-integrable kernels, are associated and hence {\bf PA} by
\cite[Chapter 1, Theorem 3.27]{Bulinski}. The finite
susceptibility condition $\sigma^2(u)\in(0,\infty)$ can be
verified by \cite[Lemma 3.71]{Karcher12}. Further examples of
random fields satisfying Theorem \ref{theo:PANA_CLT} are
\emph{stable} random fields which we consider in more detail
following \cite[Section 3.5.3]{Karcher12}.

\paragraph{\bf Max--stable random fields}

Let $X=\left\{X(t),\; t\in \Rd\right\}$ be a stationary max-stable
random field with spectral representation
$$
X(t)=\max_{i\in\N} \xi_i f_t(y_i), \quad t\in\Rd,
$$\index{max--stable random field}
where $f_t:E\to \R_+$ is a measurable function defined on the
measurable space $(E,\mu)$ for all $t\in\Rd$ with $$
\int_{E}f_t(y)\, \mu(dy)=1, \quad t\in\Rd,
$$ and $\{ (\xi_i,y_i) \}_{i\in\N}$ is the Poisson point process on
$(0,\infty)\times E$ with intensity measure $\xi^{-2}d\xi \times
\mu(dy)$. It is known that all max--stable distributions are
associated and hence {\bf PA} by \cite[Proposition
5.5.29]{Resnik08}. The field $X$ is stochastically continuous if
 $\|
f_s-f_t\|_{L^1}\to 0$ as $s\to t$ (cf. \cite[Lemma 2]{Haan84}).
Condition $\sigma^2(u)\in(0,\infty)$ is satisfied if
$$
\int_{\Rd} \int_E \min\{f_0(y), f_t(y)  \} \, \mu(dy)\, dt<\infty.
$$

\paragraph{\bf $\alpha$--stable random fields}

Let $X=\left\{X(t),\; t\in \Rd\right\}$ be a stationary
$\alpha$-stable random field ($\alpha\in (0,2)$, for simplicity
$\alpha\neq 1$) with spectral representation
$$
X(t)=\int_E f_t(x)\, \Lambda (dx), \quad t\in\Rd,
$$\index{$\alpha$--stable random field}
where $\Lambda$ is a centered independently scattered
$\alpha$--stable random measure on space $E$ with control measure
$m$ and skewness intensity $\beta: E\to [-1,1]$, $f_t:E\to \R_+$
is a measurable function on $(E,m)$ for all $t\in\Rd$. By
\cite[Proposition 3.5.1]{Samorodnitsky}, $X$ is stochastically
continuous if $\int_E |  f_s(x)-f_t(x) |^\alpha \, m(dx)\to 0$ as
$s\to t$ for any $t\in\R^d$. Condition $\sigma^2(u)\in(0,\infty)$
is satisfied if
$$
\int_{\Rd} \left(  \int_E \min\{|f_0(x)|^\alpha, |f_t(x)|^\alpha
\} \, m(dx)\right)^{1/(1+\alpha)}   \, dt<\infty.
$$

\subsection{A multivariate central limit theorem}
\label{Spo:sec:LTVol:subsec:manyfixedLevels}

If a finite number of excursion levels $u_k\in\R$, $k=1,\ldots,r$
is considered simultaneously, a multivariate analogue of Theorem
\ref{theo:methaCLT} can be proven. Introduce the notation
 $$ S_{\vec{u}}(W_n) = \left(V_d \left(A_{u_1}(X,W_n)\right),\ldots,V_d \left(A_{u_r}(X,W_n)\right)\right)^\top, \quad \Psi(\vec{u}) = \left(\Psi(u_1),\ldots,\Psi(u_r)\right)^\top .$$

\begin{theorem}[\cite{BuSpoTim12}, \cite{Karcher12}]
 Let $X$ be the above random field satisfying Theorem \ref{theo:methaCLT}. Then, for any sequence of $VH$-growing sets $W_{n}\subset\R^d $, one has
$$
      V_d \left(W_n\right)^{-1/2}\left(S_{\vec{u}}(W_n)-\Psi(\vec{u})\,V_d \left(W_n\right)\right)\stackrel{d}{\rightarrow}\mathcal{N}(0,\Sigma(\vec{u}))
$$
as $ n \rightarrow \infty $. Here,
$\Sigma(\vec{u})=(\sigma_{lm}(\vec{u}))_{l,m=1}^r$ with
\begin{equation*}
   \sigma_{lm}(\vec{u}) =\int_{\R^d} \cov\left(\ind\{X\left(0\right)\geq
u_l\},\ind\{X\left(t\right)\geq u_m\}\right)\,dt.
   \end{equation*}
\end{theorem}
If $X$ is Gaussian as in Example \ref{ex:Gauss},  we have
\begin{multline*}
   \sigma_{lm}(\vec{u})  = \\
   \frac{1}{2\pi}\int_{\mathbb{R}^d}\int_0^{\rho(t)}\frac{1}{\sqrt{1-s^2}}\exp\left\{-\frac{(u_l-a)^2-2r(u_l-a)(u_m-a)+(u_m-a)^2}{2\tau^2(1-s^2)}\right\}\,ds\,dt.\\
\end{multline*}
However, the explicit computation of the elements of matrix
$\Sigma$ for the majority of fields $X$ (except for Gaussianity)
seems to be a very complex task. In order to overcome this
difficulty in statistical applications of the methatheorem to
testing, the matrix $\Sigma$ can be (weakly) consistently
estimated from one observation of a stationary random field $X$,
see \cite{PSS10}, \cite[Section 9.8.3]{spodLNM} and references
therein.

\paragraph{\bf Statistical version of the CLT and tests}

  Let $X$ be a random field satisfying Theorem \ref{theo:methaCLT}, $u_k\in\mathbb{R}$, $k=1,\ldots,r$
  and $(W_n)_{n\in\mathbb{N}}$ be a sequence of $VH$-growing sets.
  Let $\hat{C}_n = (\hat{c}_{nlm})_{l,m=1}^r$ be a weakly consistent estimator for the nondegenerate asymptotic covariance matrix $\Sigma(\vec{u})$,
  i.e., for any $l,m=1,\ldots,r$
   $$ \hat{c}_{nlm} \stackrel{P}{\rightarrow} \sigma_{lm}(\vec{u}) \ \ \text{as } n\rightarrow\infty.$$
   Then
\begin{equation} \label{eq:StatisticalCLT}
      \hat{C}_n^{-1/2}V_d \left(W_n\right)^{-1/2}\left(S_{\vec{u}}(W_n)-\Psi(\vec{u})\,V_d \left(W_n\right)\right) \stackrel{d}{\rightarrow} \mathcal{N}(0,I).
\end{equation}
Based on the latter relation, an asymptotic test for the following
hypotheses can be constructed:
  \begin{center}
  $H_0:\; X$ \emph{is a random field satisfying Theorem \ref{theo:methaCLT} with tail distribution function } $\Psi(\cdot)$
  \end{center}\index{asymptotic test for a class of weakly dependent random fields}
vs. $H_1:$ \emph{negation of} $H_0$. As a test statistic, we use
  $$ T_n = V_d \left(W_n\right)^{-1}\left(S_{\vec{u}}(W_n)-\Psi(\vec{u})\,V_d \left(W_n\right)\right)^{\top}\hat{C}_n^{-1}\left(S_{\vec{u}}(W_n)-\Psi(\vec{u})\,V_d \left(W_n\right)\right)$$
  which is asymptotically $\chi_r^2$--distributed by continuous mapping theorem and relation \eqref{eq:StatisticalCLT}:  $T_n\xrightarrow{\mathsf{d}} \chi_r^2$ as
  $n\to\infty$. Hence, reject the null-hypothesis at a confidence level $1-\nu$ if $T_n > \chi_{r,1-\nu}^2$ where $\chi_{r,1-\nu}^2$ is the $(1-\nu)$--quantile of
  $\chi_{r}^2$--law.

\subsection{Functional limit theorems}
\label{Spo:sec:LTVol:subsec:FCLT}

A natural generalization of multivariate CLTs is a functional CLT
where the excursion level $u\in\R$ is treated as a variable, which
also appears as a (``time'') index  in the limiting Gaussian
process. In order to state the main results, introduce the
\emph{Skorokhod space} $D(\R)$ of c\`{a}dl\`{a}g functions on $\R$
endowed with the usual Skorokhod topology, cf. \cite[Section
12]{Billings99}. Denote by $\Rightarrow$ the weak convergence in
$D(\R)$.

Define the stochastic processes $Y_n=\{Y_n(u),\; u\in\R\}$ by
\begin{equation}
Y_n(u)=\frac{1}{n^{d/2}}\left( V_d\big( A_u (X, [0,n]^d) \big)
-n^d \Psi(u)  \right),\quad u\in\R.
\end{equation}
Introduce the condition
\begin{enumerate}
\item[$(\star)$] For any subset $T=\{t_1,\ldots, t_k\}\subset \Rd$
and its partition $T=T_1\cup T_2$ there exist some constants
$c(T),\gamma>0$ such that
$$
\cov\left( \prod_{t_i\in T_1} \phi_{a,b}\big(  X(t_i)\big),
\prod_{t_j\in T_2} \phi_{a,b}\big(  X(t_j)\big) \right)\le c(T)
\left(1+\dist_\infty (T_1,T_2)\right)^{-(3d+\gamma)},
$$
where $\phi_{a,b}(x)=\ind (a<x\le b) - \P (a<X(o)\le b)$ for any
real numbers $a<b$.
\end{enumerate}
 The following functional CLT is proven in
\cite[Theorem 1 and Lemma 1]{MeschenSha11}.
\begin{theorem}\label{th:FCLTVolumeMeschen}
Let $X=\{ X(t),\; t\in\Rd  \}$ be a real valued stationary random
field with a.s. continuous sample paths and a bounded density of
the distribution of $X(o)$. Let condition $(\star)$ and Theorem
\ref{theo:methaCLT} be satisfied. Then $Y_n\Rightarrow Y$ as
$n\to\infty$ where $Y=\{ Y(u), \; u\in\R\}$ is a centered Gaussian
stochastic process with covariance function
$$
C_Y(u,v)=\int_{\R^d} \cov\left(\ind\{X\left(0\right)\geq
u\},\ind\{X\left(t\right)\geq v\}\right)\,dt,\quad u,v\in\R.
$$
\end{theorem}
In particular, condition $(\star)$ is satisfied if $X\in {\bf A}$
is  square integrable with covariance function $C$ that admits a
bound
$$
|C(t)|\le \zeta\left(1+\|t\|_\infty \right)^{-\lambda}
$$
for all $t\in\Rd$ and  some $\zeta>0$, $\lambda>9d$. The proofs
 are quite technical involving a
M\'oricz bound for the moment of a supremum of (absolute values
of) partial sums of random fields on $\Z^d$, cf. \cite[Theorem
2]{Moricz83}.

For max--stable  random fields introduced in Section
\ref{Spo:sec:LTVol:subsect:fixedLevel:Stoch Cont} condition
$(\star)$ is satisfied if for any $T=\{t_1,\ldots, t_k\}\subset
\Rd$ and its partition $T=T_1\cup T_2$ there exist some constants
$c(T),\gamma>0$ such that
\begin{equation}\label{eq:CondOnF}
 \int_E \min\left\{  \max_{t_i\in  T_1} f_{t_i}(y), \max_{t_j\in  T_2} f_{t_j}(y)  \right\} \,
 \mu(dy)\le c(T)
 \left(1+\dist_\infty (T_1,T_2)\right)^{-(3d+\gamma)}.
\end{equation}

For $\alpha$--stable moving averages, i.e., $\alpha$--stable
random fields from Section
\ref{Spo:sec:LTVol:subsect:fixedLevel:Stoch Cont} with
$f_t(\cdot)=f(t-\cdot)$ for any $t\in\Rd$, condition
\eqref{eq:CondOnF} should be replaced by
\begin{multline*}
 \left( \int_{\Rd} \min\left\{  \max_{t_i\in  T_1}  f(t_i-y), \max_{t_j\in  T_2} f(t_j-y)  \right\}^\alpha \,
 m(dy) \right)^{1/(1+\alpha)}\\
 \le c(T) \left(1+\dist_\infty (T_1,T_2)\right)^{-(3d+\gamma)}.
\end{multline*}
These results  are proven (under slightly more general
assumptions) in \cite[Section 3.5.5]{Karcher12} together with
analogous conditions for infinitely divisible random fields (that
are too lengthy to give them in a review paper) as well as
examples of random fields satisfying them.

Theorem \ref{th:FCLTVolumeMeschen} together with the continuous
mapping theorem can be used to test hypotheses of Section
\ref{Spo:sec:LTVol:subsec:manyfixedLevels} with test statistic
$$
T_n=\frac{\sup_{u\in\R} Y_n(u)}{\sqrt{\E\, Y_n^2(0)}}
$$
if a large deviation result for the limiting Gaussian process $Y$
is available, cf. \cite[Corollary 1]{MeschenSha11}.

\subsection{Limit theorem for an increasing excursion level}
\label{Spo:sec:LTVol:subsec:incLevel}

If the level $u\to\infty$ one may also expect that a CLT for the
volume of the corresponding excursion set holds, provided that a
particular rate of convergence of $r$ to infinity is chosen in
accordance with the expansion rate of the observation window.

First results of this type were proven in \cite[Theorems 2.7.1,
2.7.2, 2.8.1]{IvLeon89} for stationary isotropic Gaussian random
fields with short or long range dependence. A generalization to
the case of stationary ${\bf PA}$-random fields is given in a
recent preprint \cite{DemSchmidt13}:
\begin{theorem}\label{Spo:th:VolCLTu_inf}
Let $ X=\{X(t), \; t\in\R^d\} \in {\bf PA}$ be a stationary
 random field with a continuous covariance
function $C$ such that $ |C(t)| =
\mathcal{O}\left(\left\|t\right\|_2^{-\alpha}\right)$ for some
$\alpha>3d$ as  $\left\|t\right\|_2\rightarrow\infty$. Let $X(o)$
have a bounded density $p_{X(o)}$.  Assume that the variance of
$V_d \left( A_{u_n}\big(X,[0,n]^d\big)\right)$ being equal to
\begin{equation}\nonumber
   \sigma^2_n  =\int\limits_{[0,n]^d}\int\limits_{[-x,n-x]^d}
\cov\left(\ind\{X\left(o\right)\geq
u_n\},\ind\{X\left(t\right)\geq u_n\}\right)\,dt\,
   dx
\end{equation}
satisfies\index{moving excursion level}
\begin{equation}\label{eq:CLTVolu_infVar} \sigma^2_n \to \infty,\quad n\to\infty.\end{equation}
Introduce $\gamma(x)=\sup_{y\ge x} p_{X(o)} (y)$, $x\in\R$. Choose
a sequence of excursion levels $u_n\to \infty $ such that
\begin{equation}\label{eq:CLTVolu_condVar}
\frac{n^d \gamma^{2/3}(u_n)}{\sigma^{2(\alpha+3)/3}_n} \to 0,\quad
n\to\infty.
\end{equation}
Then it holds
\begin{equation}\label{eq:CLTVolu_inf}
    \frac{V_d \left( A_{u_n}\left(X,[0,n]^d\right)\right) - n^d \P(X(o)\geq u_n)}{\sigma_n} \xrightarrow{\mathsf{d}} \mathcal{N}\left(0,1\right)
\end{equation}
as $ n \rightarrow \infty $.
\end{theorem}

Conditions \eqref{eq:CLTVolu_infVar}, \eqref{eq:CLTVolu_condVar}
are checked in \cite{DemSchmidt13} explicitly for stationary
(non-isotropic) Gaussian as well as shot noise random fields
leading to quite tractable simple expressions. For instance, it
suffices to choose $u_n=O(\sqrt{\log n})$, $n\to\infty$ in the
Gaussian case.

Student and Fisher--Snedecor random fields are considered in the
recent preprint \cite[Section 7]{LeonOl13}. CLTs for spherical
measures of excess
$$\int_{\partial B_r(o)} \ind
\{ X(t)> u(r)  \} \, {\cal H} ^{d-1}(dt)$$ of a stationary
Gaussian isotropic random field $X$ over the moving level
$u(r)\to\infty$, $r\to\infty$ are proved in \cite[Section
3.3]{Leon99}. For yet another type of geometric measures of excess
over a moving level see \cite{Leon87}.

\section{Surface area of excursion sets of stationary Gaussian random fields}
\label{Spo:sec:LTSurf}\index{surface area of excursion sets}

Limit theorems for $V_{d-1}\left(  A_u(X, W_n)\right)$ have been
first proven for one fixed level $u$ and a stationary isotropic
Gaussian random field $X$ in \cite{KratzLeon01} in dimension
$d=2$. There, the expansion of $V_{d-1}\left(  A_u(X, W_n)\right)$
in Hermite polynomials is used. In higher dimensions,  a
multivariate analogue of this result can be proven along the same
guidelines, see \cite[Proof of Theorem 1]{Shashkin13} for a
shorter proof. A CLT for the integral of a continuous function
along a level curve $\partial   A_u(X, W)$ for an a.s.
$C^1$--smooth centered mixing stationary  random field $X=\{ X(t),
\; t\in\R^2\}$ in a rectangle $W$ is proved in \cite{Iribarren89}.

\subsection{Functional limit theorems}
\label{Spo:sec:LTSurf:subsec:FCLT}

Let us focus on functional LTs for $V_{d-1}\left(  \partial A_u(X,
W_n)\right)$ proven  in \cite{MeschenSha13} for the phase space
$L^2(\R,\nu)$ (where $\nu$ is a standard Gaussian measure in $\R$)
and in \cite{Shashkin13} for the phase space $C(\R)$.

Let $X=\{ X(t),\; t\in\Rd\}$, $d>1$, be a centered stationary and
isotropic Gaussian random field with a.s. $C^1$--smooth paths and
covariance function $C\in C^2(\Rd)$ satisfying $C(o)=1$ as well as
\begin{equation}\label{eq:CondGaussXL2}
|C(t)|+ \frac{1}{1-C(t)} \sum\limits_{i=1}^d \left| \frac{\partial
C(t)}{\partial t_i} \right|+\sum\limits_{i,j=1}^d \left|
\frac{\partial^2 C(t)}{\partial t_i \partial t_j} \right|< g(t)
\end{equation}
for large $\| t\| _2$ (where $t=(t_1,\ldots, t_d)^\top$) and a
bounded continuous function \\ $g:\Rd\to\R_+$ such that $\lim_{\|
t\| _2\to \infty} g(t)=0$ and
$$
\int_{\Rd} \sqrt{g(t)}\, dt<\infty.
$$
Denote by $\nabla X (t)$ the gradient of $X(t)$. Assume that the
$(2d+2)$-dimensional random vector $\left(X(o), X(t), \nabla X
(o), \nabla X (t) \right)^\top$ is non--degenerate for all
$t\in\Rd\setminus \{ o \}$. Let $\lambda^2=-\partial^2
C(o)/\partial t_1^2$.

Introduce the sequence of random processes $\{ Y_n\}$, $n\in\N$ by
\begin{equation}\label{eq:YinL2}
Y_n(u)= \frac{2\lambda^{d/2-1}}{n^{d/2}}\left(  V_{d-1}\left(
\partial A_u(X, [0,n]^d)\right) - \E\, V_{d-1}\left(  \partial
A_u(X, [0,n]^d)\right)    \right)
\end{equation}
where $u\in\R$. They will be interpreted as random elements in
$L^2(\R,\nu)$. Let $\rightharpoonup$ denote the weak convergence
of random elements in $L^2(\R,\nu)$. Let
$$\kappa(t)=f\big( X(t) \big)
\exp\{-X^2(t)/2 \} \|  \nabla X(t)\|_2 , \quad t\in\Rd .$$

\begin{theorem}[\cite{MeschenSha13}]\label{th:FCLTSurfaceAreaL2}
Under the above assumptions on $X$ and $C$, it holds
$Y_n\rightharpoonup Y$ as $n\to\infty$ where $Y$ is a centered
Gaussian random element in $L^2(\R,\nu)$ with covariance operator
\begin{equation*}
\var \langle Y, f \rangle_{L^2(\R,\nu)} =\frac{1}{2\pi} \int_{\Rd}
\cov \left(\kappa(o),\kappa(t) \right)\, dt, \quad f\in
L^2(\R,\nu).
\end{equation*}
\end{theorem}

For $d\ge 3$, processes $Y_n$ have a continuous modification
$\tilde{Y}_n$ if conditions on $X$ starting from
\eqref{eq:CondGaussXL2} are replaced by the following ones:
\begin{enumerate}
\item Covariance function $C$ as well as all its first and second
order derivatives belong to $L^1(\R)$ \item There exist $\tau\in
(0,1)$ and $\beta>0$ such that for all $h\in [-\tau,\tau]$ and
$e_h=(h,0,0\ldots,0)^\top \in \Rd$ the determinant of the
covariance matrix of the vector
$$
\left(X(o), X(e_h),  \frac{\partial X(o)}{\partial t_1},
\frac{\partial X(e_h)}{\partial t_1} \right)^\top
$$
is not less than $|h|^\beta$.
\end{enumerate}
 Let
$\rightharpoondown$ denote the weak convergence of random elements
in $C(\R)$. Denote by $p_{X(t)}$ ($p_{X(o), X(t)}$)  the density
of $X(t)$ ( $\left(X(o), X(t)\right)^\top$), $t\in\Rd$,
respectively. Set
$$ H_t(u,v)=\E \left( \| \nabla X (o)\|_2 \| \nabla X (t)\|_2 \; |
X(o)=u,\; X(t)=v \right), \quad u,v\in\R, \quad t\in\Rd.
$$
In definition \eqref{eq:YinL2}, assume $\lambda=1$.
\begin{theorem}[\cite{Shashkin13}]\label{th:FCLTSurfaceAreaC}
Under the above assumptions on $X$ and $C$, it holds
$\tilde{Y}_n\rightharpoondown Y$ as $n\to\infty$ for $d\ge 3$
where $Y$ is a centered Gaussian random process  with covariance
function
\begin{equation*}
\cov \left(Y(u), Y(v)\right)=\int_{\Rd} \left( H_t(u,v)p_{X(o),
X(t)}(u,v)- \left(  \E \| \nabla X (o)\|_2\right)^2
p_{X(o)}(u)p_{X(t)}(v) \right)\, dt
\end{equation*}
for $u,v\in\R$.
\end{theorem}
The case $d=2$ is still open.

\section{Open problems}
\label{Spo:sec:Open}

It is a challenging problem to prove the whole spectrum of limit
theorems for $V_j\left( A_u(X,W_n) \right)$ of lower orders
$j=0,\ldots, d-2$ for isotropic $C^2$-smooth stationary Gaussian
random fields. Functional limit theorems and the case of
increasing level $u\to\infty$ are therein of special interest.
Further perspective of research is the generalisation of these
(still hypothetic) results to non--Gaussian random fields.

Another open problem is to prove limit theorems for a large class
of functionals of non--Gaussian stationary random fields that
includes the volume of excursion sets. It is quite straightforward
to do this for $$ \int_{W_n} g(X(t))\, dt$$ for a measurable
function $g:\R \to \R$ such that $\E\, g^2(X(o))<\infty$. For more
general classes of functionals of the field $X$ and the
observation window $W_n$ it is still \emph{terra incognita}.

%\begin{acknowledgement}
%If you want to include acknowledgments of assistance and the like at the end of an individual chapter please use the \verb|acknowledgement| environment -- it will automatically render Springer's preferred layout.
%\end{acknowledgement}
%
%\section*{Appendix}
%\addcontentsline{toc}{section}{Appendix}
%
%

%\input{referenc}

%----- REFERENCES -----------------------------------------------

\bibliographystyle{spmpsci}
\bibliography{abbrev,references}

\newcommand{\noopsort}[1]{} \newcommand{\printfirst}[2]{#1}
  \newcommand{\singleletter}[1]{#1} \newcommand{\switchargs}[2]{#2#1}
\begin{thebibliography}{10}
\providecommand{\url}[1]{{#1}}
\providecommand{\urlprefix}{URL }
\expandafter\ifx\csname urlstyle\endcsname\relax
  \providecommand{\doi}[1]{DOI~\discretionary{}{}{}#1}\else
  \providecommand{\doi}{DOI~\discretionary{}{}{}\begingroup
  \urlstyle{rm}\Url}\fi

\bibitem{Adler76}
Adler, R.: On generalizing the notion of upcrossings to random fields.
\newblock Adv. Appl. Probab. \textbf{8}(4), 789--805 (1976)

\bibitem{Adler81}
Adler, R.: The Geometry of Random Fields.
\newblock J. Wiley \&\ Sons, Chichester (1981)

\bibitem{AST10}
Adler, R., Samorodnitsky, G., Taylor, J.: Excursion sets of three classes of
  stable random fields.
\newblock Adv. Appl. Probab. \textbf{42}, 293--318 (2010)

\bibitem{AdlerTaylor07}
Adler, R., Taylor, J.: Random Fields and Geometry.
\newblock Springer, New York (2007)

\bibitem{ATW09}
Adler, R.J., Taylor, J.E., Worsley, K.J.: Applications of random fields and
  geometry: Foundations and case studies.
\newblock Springer Series in Statistics Springer, New York. In preparation
  (2009).
\newblock Http://webee.technion.ac.il/people/adler/hrf.pdf

\bibitem{Azais}
Azais, J.M., Wschebor, M.: Level Sets and Extrema of Random Processes and
  Fields.
\newblock Wiley, New York (2009)

\bibitem{Bel66}
Belyaev, Y.K.: On the number of intersections of a level by a {G}aussian
  stochastic process, {I}.
\newblock Theor. Probab. Appl. \textbf{11}(1), 106--113 (1966)

\bibitem{Bel67}
Belyaev, Y.K.: On the number of intersections of a level by a {G}aussian
  stochastic process, {II}.
\newblock Theor. Probab. Appl. \textbf{12}(3), 392--404 (1967)

\bibitem{Bel69}
Belyaev, Y.K., Nosko, V.P.: Characteristics of excursions above a high level
  for a {G}aussian process and its envelope.
\newblock Theor. Probab. Appl. \textbf{14}(2), 296–--309 (1969)

\bibitem{Berman92}
Berman, S.M.: Sojourns and extremes of stochastic processes.
\newblock Wadsworth \& Brooks, CA (1992)

\bibitem{Billings99}
Billingsley, P.: Convergence of Probability Measures, 2nd edn.
\newblock Wiley Series in Probability and Statistics: Probability and
  Statistics. Wiley, New York (1999)

\bibitem{BordCha12}
Bordenave, C., Chafa{\"i}, D.: Around the circular law.
\newblock Probability Surveys \textbf{9}, 1--89 (2012)

\bibitem{Bra07}
Bradley, R.C.: Introduction to {S}trong {M}ixing {C}onditions. {V}ol. 1,2,3.
\newblock Kendrick Press, Heber City, UT (2007)

\bibitem{Bul61}
Bulinskaya, E.V.: On the mean number of crossings of a level by a stationary
  {G}aussian process.
\newblock Theor. Probab. Appl. \textbf{6}, 435--438 (1961)

\bibitem{BuSpoTim12}
Bulinski, A., Spodarev, E., Timmermann, F.: Central limit theorems for the
  excursion sets volumes of weakly dependent random fields.
\newblock Bernoulli \textbf{18}, 100--118 (2012)

\bibitem{Bulinski}
Bulinski, A.V., Shashkin, A.P.: Limit theorems for associated random fields and
  related systems.
\newblock World Scientific, Singapore (2007)

\bibitem{CramerLead65}
Cram\'{e}r, H., Leadbetter, M.: The moments of the number of crossings of a
  level by a stationary normal process.
\newblock Ann. Math. Statist. \textbf{36}(6), 1656--1663 (1965)

\bibitem{Cuzick76}
Cuzick, J.: A central limit theorem for the number of zeros of a stationary
  {G}aussian process.
\newblock Ann. Probab. \textbf{4}(4), 547--556 (1976)

\bibitem{DemSchmidt13}
Demichev, V., Schmidt, J.: A central limit theorem for the volume of high
  excursions of stationary associated random fields.
\newblock Preprint (2013).
\newblock Ulm University

\bibitem{Doukhan94}
Doukhan, P.: Mixing: properties and examples, \emph{Lecture Notes in
  Statistics}, vol.~85.
\newblock Springer-Verlag, Berlin (1994)

\bibitem{Eliz88}
Elizarov, A.: Central limit theorem for the sojourn time and local time of a
  stationary process.
\newblock Theory Probab. Appl. \textbf{33}(1), 161--164 (1988)

\bibitem{fed59}
Federer, H.: Curvature measures.
\newblock Trans. Amer. Math. Soc. \textbf{93}, 418--491 (1959)

\bibitem{Haan84}
de~Haan, L.: A spectral representation for max--stable processes.
\newblock Ann. Probab. \textbf{12}(4), 1194--1204 (1984)

\bibitem{IbrZap10}
Ibragimov, I., Zaporozhets, D.: On the area of a random surface.
\newblock Zap. Nauchn. Sem. S.-Peterburg. Otdel. Mat. Inst. Steklov. (POMI)
  \textbf{312}, 154--175, arXiv:1102.3509 (2010)

\bibitem{Iribarren89}
Iribarren, I.: Asymptotic behaviour of the integral of a function on the level
  set of a mixing random field.
\newblock Probability and Mathematical Statistics \textbf{10}(1), 45--56 (1989)

\bibitem{IvLeon89}
Ivanov, A.V., Leonenko, N.N.: Statistical Analysis of Random Fields.
\newblock Kluwer, Dordrecht (1989)

\bibitem{ILRS13}
Ivanov, A.V., Leonenko, N.N., Ruiz-Medina, M.D., Savich, I.N.: Limit theorems
  for weighted nonlinear transformations of {G}aussian stationary processes
  with singular spectra.
\newblock Ann. Probab. \textbf{41}(2) (2013)

\bibitem{Kac43}
Kac, M.: On the average number of real roots of a random algebraic equation,
  vol.~49.
\newblock Bull. Amer. Math. Soc. (1943)

\bibitem{Karcher12}
Karcher, W.: On infinitely divisible random fields with an application in
  insurance.
\newblock Ph.D. thesis, Ulm University, Ulm (2012)

\bibitem{Kratz06}
Kratz, M.: Level crossings and other level functionals of stationary {G}aussian
  processes.
\newblock Probability Surveys \textbf{3}, 230--288 (2006)

\bibitem{KratzLeon01}
Kratz, M., L\'{e}on, J.: Central limit theorems for level functionals of
  stationary {G}aussian processes and fields.
\newblock J. Theor. Probab. \textbf{14}, 639--672 (2001)

\bibitem{KratzLeon10}
Kratz, M., L\'{e}on, J.: Level curves crossings and applications for {G}aussian
  models.
\newblock Extremes \textbf{13}(3), 315--351 (2010)

\bibitem{LeLiRoo83}
Leadbetter, M.R., Lindgren, G., Rootzen, H.: Extremes and Related Properties of
  Random Sequences and Processes, 1 edn.
\newblock Springer Series in Statistics. Springer (1983)

\bibitem{LeonOl13}
Leonenko, N., Olenko, A.: Sojourn measures of {S}tudent and
  {F}ischer--{S}nedecor random fields.
\newblock Bernoulli, in press (2013)

\bibitem{Leon87}
Leonenko, N.N.: Limit distributions of characteristics of exceeding a level by
  a {G}aussian field.
\newblock Mathematical notes of the Academy of Sciences of the USSR
  \textbf{41}(4), 339--345 (1987)

\bibitem{Leon88}
Leonenko, N.N.: Sharpness of the normal approximation of functionals of
  strongly correlated {G}aussian random fields.
\newblock Mathematical notes of the Academy of Sciences of the USSR
  \textbf{43}(2), 161--171 (1988)

\bibitem{Leon99}
Leonenko, N.N.: Limit theorems for random fields with singular spectrum.
\newblock Kluwer, Dordrecht (1999)

\bibitem{Liu12}
Liu, J.: Tail approximations of integrals of {G}aussian random fields.
\newblock Ann. Probab. \textbf{40}(3), 1069--1104 (2012)

\bibitem{LiuXu12}
Liu, J., Xu, G.: Some asymptotic results of {G}aussian random fields with
  varying mean functions and the associated processes.
\newblock Ann. Statist. \textbf{40}(1), 262--293 (2012)

\bibitem{Malev69}
Malevich, T.L.: Asymptotic normality of the number crossings of level zero by a
  {G}aussian process.
\newblock Theor. Probab. Appl. \textbf{14}, 287--295 (1969)

\bibitem{MarPec11}
Marinucci, D., Peccati, G.: Random Fields on the Sphere. Representation, Limit
  Theorems and Cosmological Applications, \emph{Lecture Notes of the London
  Mathematical Society}, vol. 389.
\newblock Cambridge University Press (2011)

\bibitem{MeckeStoyan02}
Mecke, K., Stoyan, D. (eds.): Morphology of Condensed Matter. Physics and
  Geometry of Spatially Complex Systems, \emph{Lecture Notes in Physics}, vol.
  600.
\newblock Springer, Berlin, Heidelberg (2002)

\bibitem{MeschenSha11}
Meschenmoser, D., Shashkin, A.: Functional central limit theorem for the volume
  of excursion sets generated by associated random fields.
\newblock Statist. Probab. Lett. \textbf{81}(6), 642–--646 (2011)

\bibitem{MeschenSha13}
Meschenmoser, D., Shashkin, A.: Functional central limit theorem for the
  measures of level surfaces of the {G}aussian random field.
\newblock Theor. Probab. Appl. \textbf{57}(1), 162--172 (2013)

\bibitem{Molch05}
Molchanov, I.: Theory of Random Sets.
\newblock Springer, London (2005)

\bibitem{Moricz83}
M{\'o}ricz, F.: A general moment inequality for the maximum of the rectangular
  partial sums of multiple series.
\newblock Acta Mathematica Hungarica \textbf{41}, 337--346 (1983)

\bibitem{PSS10}
Pantle, U., Schmidt, V., Spodarev, E.: On the estimation of integrated
  covariance functions of stationary random fields.
\newblock Scand. J. Statist. \textbf{37}, 47--66 (2010)

\bibitem{Piterbarg78}
Piterbarg, V.I.: The central limit theorem for the number of level crossings of
  a stationary {G}aussian process.
\newblock Theor. Probab. Appl. \textbf{23}(1), 178--182 (1978)

\bibitem{Resnik08}
Resnik, S.I.: Extreme values, regular variation and point processes.
\newblock Springer, Berlin, Heidelberg (2008)

\bibitem{Rice44}
Rice, S.: Mathematical analysis of random noise.
\newblock Bell System Technical Journal \textbf{23}, 282--332 (1944)

\bibitem{Rice45}
Rice, S.: Mathematical analysis of random noise.
\newblock Bell System Technical Journal \textbf{24}, 46--156 (1945)

\bibitem{Samorodnitsky}
Samorodnitsky, G., Taqqu, M.: Stable non-Gaussian random processes.
\newblock Chapman \& Hall/CRC (1994)

\bibitem{sant76}
Santal\'o, L.: Integral Geometry and Geometric Probability.
\newblock Addi\-son-Wel\-sey, Reading, Mass (1976)

\bibitem{schn93}
Schneider, R.: Convex Bodies. The {Brunn--Minkowski} Theory.
\newblock Cambridge University Press, Cambridge (1993)

\bibitem{SDT08}
Schwartzman, A., Dougherty, R.F., Taylor, J.E.: False discovery rate analysis
  of brain diffusion direction maps.
\newblock Annals of Applied Statistics \textbf{2}(1), 153--175 (2008)

\bibitem{Shashkin13}
Shashkin, A.: A functional central limit theorem for the level measure of a
  {G}aussian random field.
\newblock Statist. Probab. Lett. \textbf{83}(2), 637–--643 (2013)

\bibitem{Slud94}
Slud, E.: {MWI} representation of the number of curve--crossings by a
  differentiable {G}aussian process, with applications.
\newblock Ann. Probab. \textbf{22}(3), 1355--1380 (1994)

\bibitem{spodLNM}
Spodarev, E. (ed.): Stochastic Geometry, Spatial Statistics and Random Fields.
  Asymptotics Methods, \emph{Lecture Notes in Mathematics}, vol. 2068.
\newblock Springer, Berlin (2013)

\bibitem{SpoZap12}
Spodarev, E., Zaporozhets, D.: Asymptotics of the mean {M}inkowski functionals
  of {G}aussian excursions (2012).
\newblock Preprint

\bibitem{TayWor07}
Taylor, J.E., Worsley, K.J.: Detecting sparse signal in random fields, with an
  application to brain mapping.
\newblock Journal of the American Statistical Association \textbf{102 (479)},
  913--928 (2007)

\bibitem{Torquato02}
Torquato, S.: Random Heterogeneous Materials: Microstructure and Macroscopic
  Properties.
\newblock Springer (2002)

\bibitem{WorTay06}
Worsley, K.J., Taylor, J.: Detecting f{MRI} activation allowing for unknown
  latency of the hemodynamic response.
\newblock Neuroimage \textbf{29}, 649--654 (2006)

\end{thebibliography}

\printindex

\end{document}